\documentclass[12pt,a4paper]{article}

\usepackage{amsmath,amsfonts}
\usepackage{mathrsfs}
\usepackage{graphicx}
\def\ca#1{\ifmmode{\cal#1}\else$\cal#1$\fi}
\def \Liminf{\mathop{\underline{\lim}}\limits}

\def \zs#1{_{\lower 3pt \hbox{$\scriptstyle#1$}}}

\let\bb\mathbb       
\def\HH{{\bb H}}
\def\AA{{\bb A}}\def\UU{{\bb U}} 
\let\ca\mathcal      
\let\bb\mathbb       

\newtheorem{theorem}{Theorem}

\def\BB{{\bb B}}\def\CC{{\bb C}}

\newcommand\Pb{\mathbf{P}}
\newcommand\Ex{\mathbf{E}}

\begin{document}
\title{On Properties of   Estimators  in non Regular Situations for Poisson
Processes.} 
\author{Yury A. \textsc{Kutoyants}\\
{\small Laboratoire de Statistique et Processus, Universit\'e du Maine}}

\date{}

\maketitle

\begin{abstract}
We consider the problem of parameter estimation by observations of
 inhomogeneous Poisson process. It is well-known that if the regularity
 conditions are fulfilled then the maximum likelihood and bayesian estimators
 are consistent, asymptotically normal and asymptotically efficient. These
 regularity conditions can be roughly presented as follows: a) the intensity
 function of observed process belongs to known parametric family of functions,
 b) the model is identifiable, c) the Fisher information is positive
 continuous function, d) the intensity function is sufficiently smooth with
 respect to the unknown parameter, e) this parameter is an interior point of
 the interval. We are interested in the properties of estimators when these
 regularity conditions are not fulfilled. More precisely, we preset a review
 of the results which correspond to the rejection of these
 conditions one by one and we show how the properties of the MLE and Bayesian
 estimators change. The proofs of these results are essentially based on some
 general results by Ibragimov and Khasminskii. 

\end{abstract}

 
 \section{Introduction}

We start with the classical model of i.i.d. observations. Let $X_1, \cdots,
X_n$ be independent and identically distributed random variables with the
density function $f_*\left(x\right)$. We suppose that
$f_*\left(x\right)=f(\vartheta ,x)$, where $f\left(\cdot ,\cdot \right)$ is a 
known function depending on the unknown parameter $\vartheta\in
\Theta=\left(\alpha ,\beta \right) $. We have to estimate $\vartheta $ and to
describe the properties of estimators in the asymptotic of large samples
$\left(n\rightarrow \infty \right)$. We discuss below two estimators: maximum
likelihood and bayesian.  Let us introduce the likelihood function
$L_n\left(\vartheta ,X^n\right)=\prod_{j=1}^nf\left(\vartheta
,X_j\right)$. Then 
the maximum likelihood estimator (MLE) $\hat\vartheta _n $ and bayesian
estimator (BE) $\tilde\vartheta _n $ (for quadratic loss function and density
a priori $p\left(\cdot \right)$) are defined by the equations
\begin{equation}
\label{1}
 L_n\left(\hat\vartheta_n ,X^n\right)=\sup_{\vartheta \in
\Theta } L_n\left(\vartheta ,X^n\right),\qquad \quad  \tilde\vartheta
_n=\frac{\int_{\alpha }^{\beta } \theta p\left(\theta
\right)L_n\left(\vartheta ,X^n\right){\rm d}\theta }{\int_{\alpha }^{\beta }  p\left(\theta
\right)L_n\left(\vartheta ,X^n\right){\rm d}\theta} .
\end{equation}

It is well-known that if the {\sl conditions of  regularity} are fulfilled then
these estimators are consistent, asymptotically normal
\begin{align*}
\sqrt{n}\left(\hat\vartheta _n-\vartheta \right)\Longrightarrow {\cal
N}\left(0,{\rm I}\left(\vartheta \right)^{-1}\right),\qquad 
\sqrt{n}\left(\tilde\vartheta _n-\vartheta \right)\Longrightarrow {\cal 
N}\left(0,{\rm I}\left(\vartheta \right)^{-1}\right),
\end{align*}
and asymptotically efficient. Here  
$
{\rm I}\left(\vartheta \right)=\int_{}^{}\frac{\dot f\left(\vartheta
,x\right)^2}{f\left(\vartheta ,x\right)} \, {\rm d}\mu\left(x\right) 
$
 is the Fisher
information.  The proofs you can find in any book on asymptotical
statistics, e.g.,  Ibragimov and Khasminski (1981).

 These regularity conditions can be roughly described as follows

\begin{itemize}
\item {\sl The density $f_*\left(x\right)$ of the observed r.v.'s belongs to the
 parametric 
family, i.e., there exists a value $\vartheta _0\in \Theta=\left(\alpha ,\beta
\right) $ such that $f_*\left(x\right)=f\left(\vartheta _0,x\right)$.
\item The function $f\left(\vartheta ,x\right)$ is one or more times
differentiable w.r.t. $\vartheta $ with certain majoration of the derivatives.
\item The Fisher information ${\rm I}\left(\vartheta \right)$ is positive
 function.
\item The Fisher information ${\rm I}\left(\vartheta \right)$ is 
continuous  function.
\item The model is identifiable: if $\vartheta _1\neq \vartheta _2$ then
$f\left(\vartheta _1,x\right)\neq f\left(\vartheta _2,x\right)$.
\item The true value $\vartheta _0$ is an interior point of the set $\Theta
$, i.e., $\vartheta _0\neq \alpha $ and $\vartheta _0\neq \beta $.
\item We can observe all values the random variables $X_j$.
\item The statistical  model is fixed and can not be chosen in some optimal
way. }
\end{itemize}
Of course, this list is not exhaustive and the other conditions can be
 mentioned too.  We are interested by the properties of estimators, when the
 similar regularity conditions are not fulfilled for some models of continuous
 time stochastic processes. More  precisely, we replace these regularity
 conditions by other conditions and study the properties of estimators under
 these new conditions. This approach allows to understand better the role
 of each regularity condition in the properties of estimators. As the model of
 observations in this work we take inhomogeneous Poisson process. The similar
 work concerning parameter estimation for ergodic diffusion processes was
 already published (see \cite{K08}), but it seems that the more detailed
 exposition of the proofs will be useful and it is given here.

\section{Regular case}

We observe $n$ independent trajectories $X^n=\left(X_1,\ldots,X_n\right)$,
 where $X_j=$ $\{X_j\left(t\right),$ $ 0\leq t\leq \tau\}$, of a Poisson
 process $X^\tau =\left\{X\left(t\right), 0\leq t\leq \tau \right\}$ of intensity
 function $\lambda _*=\left\{\lambda _*\left(t\right), 0\leq
 t\leq \tau \right\}$, i.e., $X\left(0\right)=0$, the increments on disjoint
 intervals are independent and
$$
\Pb\left\{ X\left(t\right)=k\right\} =\frac{\Lambda
\left(t\right)^k}{k!}\exp\left\{-\Lambda\left(t\right) \right\},
\qquad \quad \Lambda\left(t\right)= \int_{0}^{t}\lambda_*
\left(s\right)\;{\rm d}s. 
$$
The Poisson process
 sometimes is   defined as a series of events
 $0<t_1<t_2<\ldots <t_M<T$ and $X\left(t\right), 0\leq t\leq \tau $  is the
 corresponding counting process, i.e., $X\left(t\right)$ is equal to the
 number of events observed up to time $t$. The process $X^\tau $ is c\`adl\`ag (right
 continuous with left limits at every point $t$). 

The same model of observation
 we obtain in the case of $\tau $-periodic Poisson process
 $X^{T_n}=\left\{X\left(s\right), 0\leq s\leq T_n\right\}$, if the intensity
 function $\lambda _*\left(s\right)$ is $\tau $-periodic and $T_n=\tau
 n$. Then we can cut the trajectory $X^{T_n}$ on $n$ pieces
 $X_j\left(t\right)=X\left(t+\left(j-1\right)\tau
 \right)-X\left(\left(j-1\right)\tau \right), 0\leq t\leq \tau  $  with
 $j=1,\ldots,n$. If we suppose that the period $\tau $ is known (does not
 depend on $\vartheta $). As the increments of the Poisson process are independent,
 this model coincides with the mentioned above one.

 The 
 Statistician can suppose that this intensity function  belongs to some
parametric class of functions, 
i.e., $\lambda_* =\lambda_\vartheta  $, where   $\lambda_\vartheta 
=\left\{\lambda \left(\vartheta ,t\right),\right.$ $\left. 0\leq t\leq 
\tau\right\}$ with 
$\vartheta\in \Theta =\left(\alpha ,\beta \right)$. 
Therefore he (or she) obtains the problem of estimation of the parameter $\vartheta $ by
the observations $X^n$ of the Poisson process of intensity function $\lambda
_{\vartheta }, \vartheta \in \Theta $.

We suppose that the intensity  is bounded positive function and hence the
likelihood ratio function for this parametric family is 
\begin{equation}
\label{LR}
L\left(\vartheta ,X^n\right)=\exp\left\{\sum_{j=1}^{n}\int_{0}^{\tau }\ln
\lambda \left(\vartheta ,t\right)\,{\rm d}X_j\left(t\right)-n \int_{0}^{\tau
}\left[\lambda \left(\vartheta ,t\right)-1 \right]\,{\rm d}t\right\}
\end{equation}
and the MLE $\hat \vartheta _n$ and BE $\tilde \vartheta _n$ for quadratic
loss function and prior 
density $p\left(\theta \right),\theta \in \Theta $ (positive, continuous on
$\Theta $) are defined by the same equations \eqref{1}.

 {\bf Regularity Conditions:}
\begin{enumerate}
\item {\sl There exists $\vartheta _0\in
\Theta $ such that  $\lambda _*\left(t\right)=\lambda\left(\vartheta_0
,t\right),0\leq t\leq \tau $.} 
\item {\sl The function $\sqrt{\lambda \left(\vartheta ,t\right)},0\leq t\leq
\tau $ has two continuous bounded derivatives with respect to $\vartheta $.
\item {\sl The Fisher information}
$$
0<{\rm I}\left(\vartheta \right)=\int_{0}^{\tau }\frac{\dot \lambda
\left(\vartheta ,t\right)^2}{\lambda \left(\vartheta ,t\right)}\;{\rm d}t <\infty . 
$$
\item  The Fisher information
${\rm I}\left(\vartheta \right)$ is continuous function.}
\item {\sl The condition of identifiability is fulfilled: for any $\nu >0$}
\begin{equation*}
\label{id}
\inf_{\left|\theta -\vartheta _0\right|>\nu }\int_{0}^{\tau
}\left[\sqrt{\lambda \left(\vartheta ,t\right)}-\sqrt{\lambda
\left(\vartheta_0 ,t\right)} \right]^2\;{\rm d}t>0 . 
\end{equation*}
\item {\sl The parameter $\vartheta _0$ is an interior point of the set
$\Theta =\left(\alpha ,\beta \right)$.} 
\item {\sl The process $X_j\left(t\right)$ is observed on the whole interval
$\left[0,\tau \right]$. 
\item  The model of observed process is fixed, i.e., in the statement of the
problem the intensity function $\lambda _\vartheta $ is given (can not be
chosen by the statistician).}
\end{enumerate}
 
Of course, Condition 2 implies 4, but we present both of them, because we
consider below the case, when 4 is not fulfilled. The properties of estimators
are described in the following theorem.

\begin{theorem}
\label{T1} Let the conditions of regularity be fulfilled, then the MLE
$\hat\vartheta _n$ and the BE $\tilde\vartheta _n$
are consistent, asymptotically normal
\begin{equation*}
\sqrt{n}\left(\hat\vartheta _n-\vartheta_0 \right)\Longrightarrow {\cal
N}\left(0,\frac{1}{{\rm I}\left(\vartheta_0 \right)} \right),
\qquad 
\sqrt{n}\left(\tilde\vartheta _n-\vartheta_0 \right)\Longrightarrow {\cal 
N}\left(0,\frac{1}{{\rm I}\left(\vartheta_0 \right)} \right)
\end{equation*}
  asymptotically efficient  and  the moments of these estimators converge too.
\end{theorem}
The proof can be found in \cite{Kut98}, Theorems 2.4 and 2.5. 

  This proof is essentialy based on the general results obtained by Ibragimov
and Khasminskii \cite{IH81}, which we present below in a bit more general
situation, than we need for this theorem.  Let us denote by $Z_n\left(u\right)$
the normalized likelihood ratio process
$$
Z_n\left(u\right)=\frac{L\left(\vartheta
_0+\varphi _n u,X^n\right)}{L\left(\vartheta _0,X^n\right)}
, \qquad u\in \UU_n=\left(\frac{\left(\alpha -\vartheta
_0\right)}{\varphi _n},\frac{\left(\beta  -\vartheta _0\right)}{\varphi _n}\right)
$$
where $\varphi _n\rightarrow 0$ and the rate of this convergence is such that
$Z_n\left(u\right)$ has some non degenerate limit (in distribution)
$Z\left(u\right)$.  Below we suppose that in the bayesian case  ($\vartheta $ is a
random variable) the loss
function is quadratic and the density {\sl  a priory}
$p\left(\vartheta \right), \vartheta \in \left(\alpha ,\beta\right) $ is
continuous positive function. Let us  define the random variables $\hat u$ and
$\tilde u$ by the relations
\begin{equation}
\label{def}
Z\left(\hat u\right)=\sup_{u\in\cal R} Z\left(u\right),\qquad \quad \tilde
u=\frac{\int_{\cal R}^{}u\,Z\left(u\right)\,{\rm
d}u}{\int_{\cal R}^{}u\,Z\left(u\right)\,{\rm d}u} .
\end{equation}
The study of the likelihood ratio $Z_n\left(\cdot \right)$
allows to describe the properties of estimators (maximum likelihood and
bayesian) and this is illustrated by the following theorem.
\begin{theorem}
\label{T2} {\rm (Ibragimov, Khasminskii)}
 Suppose that the following conditions are fulfilled
\begin{enumerate}
\item There exist   constants $a>1, B>0$, such that for all $u\in \UU_n$
\begin{equation}
\label{A}
\Ex_\vartheta
\left|Z_n^{1/2}\left(u_2\right)-Z_n^{1/2}\left(u_1\right) \right|^{2}\leq
B\,\left|u_2-u_1\right|^{a}
\end{equation}
\item There exist constants $\kappa >0$ and $\gamma >0$ such that for all $u\in \UU_n$
\begin{equation}
\label{B}
\Ex_\vartheta Z_n^{1/2}\left(u\right)\leq e^{-\kappa \left|u\right|^\gamma }
\end{equation}
\item
 The marginal distributions 
$$
\left(Z_n\left(u_1\right),\ldots,Z_n\left(u_k\right)\right)\Longrightarrow
\left(Z\left(u_1\right),\ldots,Z\left(u_k\right)\right)
$$
 and $Z\left(\cdot
\right)$ attains  with probability 1 its maximal value at a 
unique point $\hat u$.
\end{enumerate}
Then, the MLE $\hat\vartheta _n$ and BE $\tilde\vartheta _n$ are consistent,
$$
\varphi _n^{-1}\left(\hat \vartheta _n-\vartheta \right)\Longrightarrow \hat
u,\qquad 
\varphi _n^{-1}\left(\tilde \vartheta _n-\vartheta \right)\Longrightarrow \tilde
u ,
 $$
and for any $p>0$
$$
\Ex_\vartheta \left|\frac{\hat \vartheta 
_n-\vartheta }{\varphi _n}\right|^p \longrightarrow \Ex_\vartheta \left|\hat
u\right|^p,\qquad \Ex_\vartheta \left|\frac{\tilde \vartheta 
_n-\vartheta }{\varphi _n}\right|^p \longrightarrow \Ex_\vartheta \left|\tilde
u\right|^p .
$$
\end{theorem}
For the proof (essentially more general results) see  \cite{IH81}, Theorems
3.1.1 and 3.2.1. Note that in the case of bayesian estimators it is sufficient that the
parameter $a>0$. 

\bigskip
 In the regular case of the Theorem \ref{T1} the sequence $\varphi
_n=n^{-1/2}$ and the limit process is
$$
  Z\left(u\right) =\exp\left\{u\;\zeta \left(\vartheta
_0\right)-\frac{u^2}{2}\;{\rm I}\left(\vartheta_0 \right)\right\},\qquad u\in
{\cal R},
$$
where $\zeta \left(\vartheta _0\right)\sim {\cal N}\left(0, {\rm
I}\left(\vartheta_0 \right)\right)$.  Hence 
\begin{equation*}
\hat u=\frac{\zeta \left(\vartheta _0\right)}{{\rm
I}\left(\vartheta_0 \right)}\quad \sim \quad {\cal N}\left(0, {\rm
I}\left(\vartheta_0 \right)^{-1}\right)
\end{equation*}

To check the conditions \eqref{A} and \eqref{B} in the case of inhomogeneous
Poisson processes we use the following estimates (below $\vartheta _0$ is the
true value and  $\vartheta _i=\vartheta
_0+\frac{u_i}{\sqrt{n}}$)
\begin{align}
\label{C}
&\Ex_{\vartheta_0}
\left|Z_n^{1/2}\left(u_2\right)-Z_n^{1/2}\left(u_1\right) \right|^{2}=2-2\Ex_{\vartheta_0}
\left[Z_n\left(u_2\right)Z_n\left(u_1\right)\right]^{1/2}\nonumber\\
&\qquad \qquad  =
2-2\Ex_{\vartheta_1}\left[\frac{Z_n\left(u_2\right)}{Z_n\left(u_1\right)}
\right]^{1/2}  \nonumber \\
&\qquad \qquad =
2-2\exp\left\{-\frac{n}{2}\int_{0}^{\tau }\left[\sqrt{\lambda \left(\vartheta
_2,t\right)} -\sqrt{\lambda \left(\vartheta _1,t\right)}\right]^{2}{\rm
d}t\right\}\nonumber\\ 
&\qquad \qquad  \leq n\int_{0}^{\tau }\left[\sqrt{\lambda \left(\vartheta
_2,t\right)} -\sqrt{\lambda \left(\vartheta _1,t\right)}\right]^{2}{\rm
d}t
\end{align}
and $\left( \vartheta _u=\vartheta_0+\frac{u}{\sqrt{n}}\right)$
\begin{align}
\label{D}
&\Ex_{\vartheta_0} Z_n^{1/2}\left(u\right)\nonumber\\
& =\left(\Ex_{\vartheta_0} \exp\left\{
\frac{1}{2}\int_{0}^{\tau } \ln\frac{\lambda \left(\vartheta
_u,t\right)}{\lambda \left(\vartheta_0,t\right)}{\rm
d}X\left(t\right)-\frac{1}{2}\int_{0}^{\tau }\left[\lambda
\left(\vartheta_u,t\right)-\lambda \left(\vartheta_0,t\right)  \right]{\rm d}t
\right\}\right)^n\nonumber\\
& =
\exp\left\{-\frac{n}{2}\int_{0}^{\tau }\left[\sqrt{\lambda \left(\vartheta
_u,t\right)} -\sqrt{\lambda \left(\vartheta _0,t\right)}\right]^{2}{\rm
d}t\right\}.
\end{align}
The regularity conditions allow to obtain the low and upper estimates
\begin{align}
\label{E}
c\left|\vartheta _2-\vartheta _1\right|^2\leq \int_{0}^{\tau
}\left[\sqrt{\lambda \left(\vartheta 
_2,t\right)} -\sqrt{\lambda \left(\vartheta _1,t\right)}\right]^{2}{\rm
d}t\leq C\left|\vartheta _2-\vartheta _1\right|^2
\end{align}
which provide immediately \eqref{A} and \eqref{B}. Using the direct expansion of the
functions 
$$
\lambda \left(\vartheta _0+\frac{u}{\sqrt{n}},t\right)=\lambda
\left(\vartheta _0,t\right)+\frac{u}{\sqrt{n}}\dot \lambda \left(\vartheta
_0,t\right)+o\left(\frac{u}{\sqrt{n}} \right)
 $$
and $ \ln \lambda \left(\vartheta _0+\frac{u}{\sqrt{n}},t\right)$ we obtain
the following representation of the likelihood ratio 
\begin{align*}
Z_n\left(u\right)=\exp\left\{u\Delta _n\left(\vartheta
_0,X^n\right)-\frac{u^2}{2}{\rm I}\left(\vartheta _0\right)+r_n\right\},
\end{align*}
where
$$
\Delta_n\left(\vartheta_0,X^n\right)=\frac{1}{\sqrt{n}}\sum_{j=1}^{n}\int_{0}^{\tau
}\frac{\dot \lambda \left(\vartheta_0,t\right) }{ \lambda \left(\vartheta
_0,t\right)}\left[{\rm d}X_j\left(t\right)-\lambda \left(\vartheta
_0,t\right){\rm d}t\right]\Longrightarrow \zeta \left(\vartheta_0\right)
$$
and $r_n\rightarrow 0$. This representation provides the convergence of the
marginal distributions of the process $Z_n\left(\cdot \right)$ to the marginal
distributions of the process $Z\left(\cdot \right)$.  Therefore all conditions
of the Theorem \ref{T2} are fulfilled and the MLE and BE are consistent,
asymptotically normal. Let us remind how the weak convergence of the
likelihood ratio process provides these properties of estimators.

Suppose that we already have the weak
convergence of the stochastic processes
\begin{equation}
\label{6}
Z_n\left(\cdot \right) \Longrightarrow Z\left(\cdot \right)
\end{equation}
 in the space
of continuous on ${\cal R}$ functions vanishing in infinity.
 Then according to \cite{IH81} the asymptotic
normality of the MLE can be obtained by the following way.
\begin{align}
&\Pb\left\{\sqrt{n}\left(\hat\vartheta _n-\vartheta_0
\right)<x \right\}=\Pb\left\{\sup_{\sqrt{n}\left(\theta
-\vartheta _0\right)<x}L\left(\vartheta,X^n\right) >\sup_{\sqrt{n}\left(\theta
-\vartheta _0\right)\geq x}L\left(\vartheta,X^n\right) \right\}\nonumber\\ 
&\quad=\Pb\left\{\sup_{\sqrt{n}\left(\theta
-\vartheta
_0\right)<x}\frac{L\left(\vartheta,X^n\right)}{L\left(\vartheta_0,X^n\right) }
>\sup_{\sqrt{n}\left(\theta 
-\vartheta _0\right)\geq
x}\frac{L\left(\vartheta,X^n\right)}{L\left(\vartheta_0,X^n\right) }
\right\}\nonumber\\ 
&\quad=\Pb\left\{\sup_{u<x}Z_n\left(u\right) >\sup_{u\geq
x}Z_n\left(u\right) \right\}\longrightarrow
\Pb\left\{\sup_{u<x}Z\left(u\right) >\sup_{u\geq x}Z\left(u\right)
\right\}\nonumber\\ 
&\quad =\Pb\left(\frac{\zeta \left(\vartheta _0\right) }{{\rm
I}\left(\vartheta_0 \right)}<x\right),\quad {\rm i.e.}\quad
\sqrt{n}\left(\hat\vartheta _n-\vartheta_0 \right)\Longrightarrow {\cal
N}\left(0,\frac{1}{ {\rm I}\left(\vartheta_0 \right)}\right).
\label{8}
\end{align}
where we put $\vartheta =\vartheta_0+u/\sqrt{n}$.

For the BE we  change the variable $\theta =\vartheta
_0+u/\sqrt{n}\equiv \vartheta _u$
\begin{eqnarray*}
&& \tilde\vartheta _n=\frac{\int_{\alpha }^{\beta }\theta p\left(\theta
\right)L\left(\theta ,X^n\right){\rm d}\theta }{\int_{\alpha }^{\beta } p\left(\theta
\right)L\left(\theta ,X^n\right){\rm d}\theta } =\vartheta
_0+\frac{1}{\sqrt{n}} \frac{\int_{\UU_n }u p\left(\vartheta _u
\right)L\left(\vartheta _u ,X^n\right){\rm d}u }{\int_{\UU_n } p\left(\vartheta _u
\right)L\left(\vartheta _u ,X^n\right){\rm d}u },
\end{eqnarray*}
 Then using  the convergence $p\left(\vartheta
_u \right)\rightarrow 
p\left(\vartheta _0\right)$  (according to \cite{IH81})  we can write
\begin{align}
&\Pb_{\vartheta _0}\left\{\sqrt{n}\left(\tilde\vartheta _n-\vartheta_0
\right)<x \right\}=\Pb\left\{ \frac{\int_{\UU_n }u\; p\left(\vartheta _u
\right)Z_n\left(u\right)\;{\rm d}u }{\int_{\UU_n } p\left(\vartheta _u
\right)Z_n\left(u\right)\;{\rm d}u }<x \right\}\nonumber\\ 
&\qquad
\longrightarrow \Pb\left\{ \frac{\int_{R}u\; Z\left(u \right)\;{\rm d}u
}{\int_{R} Z\left(u \right)\;{\rm d}u }<x \right\} =\Pb\left(\frac{\zeta
\left(\vartheta _0\right) }{{\rm I}\left(\vartheta_0 \right)}<x\right)
\label{9}
\end{align}
because the elementary calculus yield the equality
$$
\int_{R}u \;Z\left(u \right)\;{\rm d}u=\int_{R}^{}u\,e^{u\zeta \left(\vartheta
_0\right)-\frac{u^2}{2 }{\rm I}\left(\vartheta_0
\right)}\;{\rm d}u=\frac{\zeta \left(\vartheta _0\right) }{{\rm I}\left(\vartheta_0
\right)}\;\int_{R}Z\left(u \right)\;{\rm d}u.
$$
 
Hence 
$$
\sqrt{n}\left(\tilde\vartheta _n-\vartheta_0\right)\Longrightarrow {\cal
N}\left(0, \frac{1}{{\rm I}\left(\vartheta_0 \right)}\right).
$$

Moreover, by Theorem \ref{T2}
$$
\left(n{\rm I}\left(\vartheta_0 \right)\right) ^{\frac{p}{2}}\Ex_{\vartheta_0
}\left|\hat\vartheta _n-\vartheta_0 
\right|^p\longrightarrow   \Ex\left|\zeta \right|^p, \quad \left(n{\rm
I}\left(\vartheta_0 \right)\right) ^{\frac{p}{2}}\Ex_{\vartheta_0 
}\left|\tilde\vartheta _n-\vartheta_0 
\right|^p\longrightarrow   \Ex\left|\zeta  \right|^p ,
$$
 where $\zeta \sim {\cal N}\left(0,1\right)$

\section{Misspecified model}

Suppose now that the parametric family $\left\{\lambda_\vartheta
,\vartheta \in \Theta \right\}$ does not
correspond to the observed process $X^n$, i.e., the value $\vartheta _0\in
\Theta $, such that $\lambda _*=\lambda _{\vartheta _0}$ does not exist,
 but the statistician nevertheless uses this model to estimate the
parameter $\vartheta $ ({\sl no true model} case), i.e., he (or she) calculates
the likelihood ratio function by \eqref{LR}, where  $X^n$ are  observations of
the Poisson process of intensity function $\lambda _*\left(\cdot \right)$.
It can be shown that the MLE and BE converge to the value
\begin{equation}
\label{KL}
\vartheta _*=\arg\inf_{\theta\in \Theta  }\int_{0}^{\tau
}\left[\frac{\lambda \left(\vartheta ,t\right)}{\lambda_*\left(t\right) }-1-
\ln \frac{\lambda\left(\vartheta
,t\right)}{\lambda_*\left(t\right)}\right]\lambda_*\left(t\right)\;{\rm d}t, 
\end{equation}
which minimizes the Kullback-Liebler distance between the measure $\Pb_*$,
which corresponds to the observed process with intensity $\lambda _*$
and the parametric  family  $\left\{\Pb_\vartheta ,\vartheta \in\Theta
\right\}$.
Note that if $\lambda _*\left(t\right)=\lambda \left(\vartheta
 _0,t\right),$ $ 0\leq t\leq \tau  $, then
  $\vartheta _*=\vartheta _0 $, i.e., the { both estimators are
 consistent}.

 Moreover if $\vartheta _* $ is an interior point of the set $\Theta $, then
 these estimators are asymptotically normal:
$$
\sqrt{n}\left(\hat\vartheta _n-\vartheta _*\right)\Longrightarrow {\cal
N}\left(0,D_*^2\right),\qquad  \sqrt{n}\left(\tilde\vartheta _n-\vartheta
_*\right)\Longrightarrow {\cal N}\left(0,D_*^2\right).
$$
 Here $D_*^2=d_*^2\;{\rm
I}_*^{-2} $ with
$$
d_*^2=\int_{0}^{\tau }\frac{\dot \lambda \left(\vartheta
_*,t\right)^2}{\lambda \left(\vartheta _*,t\right)^2}\lambda_*
\left(t\right)\,{\rm d}t,\quad  {\rm I}_*=d_*^2+\int_{0}^{\tau }\ddot \lambda \left(\vartheta
_*,t\right)\left[1-\frac{\lambda_*
\left(t\right)}{\lambda \left(\vartheta _*,t\right)}\right]{\rm d}t . 
$$
Note that in this case  the pseudo-LR function $Z_n\left(u\right)$ constructed
on the base of the wrong parametric model has a different limit
$$
Z_n\left(u\right)=\frac{L\left(\vartheta
_*+u/\sqrt{n},X^T\right)}{L\left(\vartheta _*,X^T\right) }\Longrightarrow
Z\left(u\right)=\exp\left\{u\,\zeta _*-\frac{u^2}{2}\,{\rm I}_*\right\} 
$$
where $\zeta _*\sim {\cal N}\left(0,d_*^2\right) $.
The details of this proof can be found in \cite{Kut98}. See as well  Yoshida
and Hayashi \cite{YH}. 

We are interested here by a different problem. The intensity of observed
 process $ \lambda_* \left(t\right)$ can be written as contaminated version of
 the parametric model $\lambda_* \left(t\right)=\lambda \left(\vartheta
 _0,t\right)+h\left(t\right),\quad 0\leq t\leq T,$ where $h\left(\cdot
 \right)$ (contamination) is unknown function. Hence $\vartheta _*=\vartheta
 _*\left(h\right)$ is the point of the minimum of the Kullback-Leibler
 distance \eqref{KL}.  We can put the following question: 
\begin{center}
{\bf when $\vartheta
 _*=\vartheta _0$} ?
\end{center}
 i.e., when nevertheless the MLE and BE are consistent?

 We consider two
 situations. The first one ({\bf smooth}), when the support $\AA\subset
 \left[0,\tau \right]$ 
 of the function $h\left(\cdot \right)$ is known and $\AA^c=\left[0,\tau
 \right]\setminus \AA\neq \emptyset
 $. We can modify the likelihood ratio and  write it as 
\begin{align*}
\ln L\left(\vartheta ,X^n\right)&=\sum_{j=1}^{n}\int_{0}^{\tau }\ln \lambda
\left(\vartheta ,t\right)\;1_{\left\{t \in \AA^c\right\}}{\rm
d}X_j\left(t\right) -n \int_{0}^{\tau }\left[\lambda \left(\vartheta
,t\right)-1\right]1_{\left\{t \in \AA^c\right\}} {\rm d}t ,
\end{align*}
i.e., we exclude the observations on $\AA $ and  define the MLE $\hat\vartheta
 _n$ and BE $\tilde\vartheta _n$ with the 
 help of this function (we call them  
 pseudo-MLE and pseudo-BE). Then we have to check if the set of intensity  functions
$\left\{\lambda \left(\vartheta ,t\right), t\in \AA^c, \vartheta \in \Theta \right\}$
satisfies the  correspondingly modified regularity conditions. For example,
 the Fisher information 
$$ 
{\rm I}_*\left(\vartheta \right)=\int_{\AA^c}^{}
 \frac{\dot \lambda \left(\vartheta ,t\right)^2}{\lambda \left(\vartheta
 ,t\right)}\;{\rm d}t >0
$$ 
and  the condition of
 identifiability : for any $\nu >0$
\begin{equation*}
\inf_{\left|\theta -\vartheta _0\right|>\nu }\int_{\AA^c}^{
}\left[\sqrt{\lambda \left(\vartheta ,t\right)}-\sqrt{\lambda
\left(\vartheta_0 ,t\right)} \right]^2\;{\rm d}t>0  
\end{equation*}
If these conditions  are fulfilled, then the estimators $\hat\vartheta _n$ and
 $\tilde\vartheta _n$  converge to the true value (are consistent) and are
 asymptotically normal.

{\bf Discontinuous intensity functions.}  Suppose that intensity of the
observed process is
\begin{equation*}
\lambda
_*\left(t\right)=\left[g_1\left(t\right)+h_1\left(t\right)\right]\,
1_{\left\{t<\vartheta_0\right\}
}+\left[g_2\left(t\right)+h_2\left(t\right)\right]\,1_{\left\{t\geq \vartheta_0
\right\}},   
\end{equation*}
where $g_1\left(\cdot \right)<g_2\left(\cdot \right)$ are known positive
functions and the functions  $h_1\left(\cdot
\right),$ $h_2\left(\cdot \right)$ are unknown.  We have to estimate the time
$\vartheta _0$ of switching of intensity function (change point estimation problem).
The MLE and BE are constructed on the base of the model  with
\begin{equation*}
\lambda \left(\vartheta ,t\right)=g_1\left(t\right)\,1_{\left\{t<\vartheta\right\}
}+g_2\left(t\right)\,1_{\left\{t\geq \vartheta
\right\}},\qquad 0\leq t\leq 
\tau,   
\end{equation*}
with the likelihood ratio function \eqref{LR}, i.e. as if 
$h_i\left(t\right)\equiv 0$, but the  observations $X^n$ used in \eqref{LR} 
contain, of course, $h_i\left(\cdot \right)$. The Kullback-Leibler distance \eqref{KL}
for $\vartheta <\vartheta _0$ is
\begin{align*}
J_{KL}\left(\vartheta \right)&=\int_{0}^{\vartheta
}\left[\frac{g_1\left(t\right)}{g_1\left(t\right)+h_1\left(t\right)}-1-\ln
\frac{g_1\left(t\right)}{g_1\left(t\right)+h_1\left(t\right)}\right]
\left[g_1\left(t\right)+h_1\left(t\right)\right]{\rm d}t \\
&\quad + \int_{\vartheta }^{\vartheta_0
}\left[\frac{g_1\left(t\right)}{g_2\left(t\right)+h_1\left(t\right)}-1-\ln
\frac{g_2\left(t\right)}{g_1\left(t\right)+h_1\left(t\right)}\right]
\left[g_1\left(t\right)+h_1\left(t\right)\right]{\rm d}t \\
&\quad + \int_{\vartheta_0 }^{\tau 
}\left[\frac{g_2\left(t\right)}{g_2\left(t\right)+h_2\left(t\right)}-1-\ln
\frac{g_2\left(t\right)}{g_2\left(t\right)+h_2\left(t\right)}\right]
\left[g_2\left(t\right)+h_2\left(t\right)\right]{\rm d}t
\end{align*}
and the similar expression we have for $\vartheta >\vartheta _0$. It is easy to see
that if the functions
$h_i\left(\cdot \right)$ satisfy the following condition 
\begin{equation}
\label{13}
0<g_1\left(t \right)+h_1\left(t\right)<\frac{g_2\left(t \right)-g_1\left(t
\right)}{\ln\frac{ g_2\left(t \right)}{g_1\left(t \right)}}< g_2\left(t
\right)+h_2\left(t\right),
\end{equation}  
then 
$$
\left.\frac{{\rm d} J_{KL}\left(\vartheta \right)}{{\rm d}\vartheta
}\right|_{\vartheta <\vartheta _0} <0,\quad {\rm and}\quad \left.\frac{{\rm
d}J_{KL}\left(\vartheta \right) }{{\rm d}\vartheta 
}\right|_{\vartheta >\vartheta _0} >0.
$$
Hence the minimum of this function is reached at the point $\vartheta
_*=\vartheta _0$ and this provides the consistency of the estimators $\hat
\vartheta _n$ and $\tilde \vartheta _n$. 
If we denote
$x=g_2\left(t\right)/g_1\left(t\right)$,
$h_i=h_i\left(t\right)/g_1\left(t\right)$, then we obtain the following
{\sl regions of consistency} for $h_i$
$$
h_1<\frac{x-1}{\ln x}-1,\qquad h_2>\frac{x-1}{\ln x}-x.
$$ 
It is important to note that the values of $h_i$ can be sufficiently large.

It can be shown that the rate of
convergence is essentially better than in regular case, and
$n\left(\hat\vartheta _n-\vartheta _0\right) $ converges in distribution to
some random variable (see similar results in Dabye and Kutoyants \cite{DK},
and Dabye, Farinetto and Kutoyants \cite{DFK}).

\section{Non identifiable model}

Suppose that we have the same model for the different values of the parameter,
i.e., $\lambda \left(\vartheta _1,t\right)=\lambda \left(\vartheta _l,t\right),
l=2,\ldots,k$, where $\vartheta _l\neq \vartheta _i, l\neq i$ and $\vartheta
_l,\vartheta _i\in \Theta $ (too many true models). It is well-known that the
MLE converges to the set $\left\{\vartheta _1,\ldots,\vartheta _k\right\}$ of
all {\sl true values}.

Let us introduce the Gaussian vector ${\bf \zeta}=\left(\zeta
_1,\ldots,\zeta _k\right) $ with zero mean and covariance matrix $\varrho
=\left(\varrho _{li}\right)$ 
$$
\varrho _{li}=\Ex\left(\zeta _l\zeta _i\right)=\left({\rm I}\left(\vartheta
_l\right){\rm I}\left(\vartheta _i\right)\right)^{-1/2}
\int_{0}^{\tau }\frac{\dot \lambda \left(\vartheta _l,t \right)\dot
\lambda \left(\vartheta _i,t \right)}{\lambda  \left(\vartheta _i,t \right)} {\rm d}t
$$
where the Fisher informations 
$$
{\rm I}\left(\vartheta_l\right)=\int_{0}^{\tau }\frac{\dot
\lambda \left(\vartheta _l,t \right)^2}{\lambda \left(\vartheta _l,t\right)}{\rm d}t>0, \qquad
l=1,2,\ldots,k .
$$

Define two random variables: discrete and continuous
$
\hat\vartheta =\sum_{l=1}^{k}\vartheta _l\,1_{\left\{\HH_l\right\}},$ and $
\tilde\vartheta =\sum_{l=1}^{k}\vartheta _l\,Q_l ,
$
where (we suppose that ${\Pb}\left\{\left|\zeta _l\right|=\left|\zeta
_i\right|\right\}=0 $)
$$
\HH_l=\left\{\omega :\left|\zeta _l\right|>\max_{i\neq l}\left|\zeta
_i\right|\right\},\qquad  Q_l=\frac{p\left(\vartheta _l\right){\rm
I}\left(\vartheta _l\right)^{-1/2} e^{\zeta
_l^2/2}}{\sum_{i=1}^{k}p\left(\vartheta _i\right){\rm 
I}\left(\vartheta _i\right)^{-1/2} e^{\zeta _l^2/2}}.
$$
It can be shown  that the MLE and BE have the following limits: $$\hat\vartheta
_n \Longrightarrow  \hat\vartheta,\qquad \qquad  \tilde\vartheta _n\Longrightarrow
\tilde\vartheta.
 $$
 Moreover 
$$
\sqrt{n}\left(\hat\vartheta_n - \hat\theta_n\right)\Longrightarrow \hat\zeta
,\qquad \quad  \sqrt{n}\left(\tilde\vartheta_n -
\tilde\theta_n\right)\Longrightarrow \tilde\zeta, 
$$
where $\hat\theta_n,\tilde\theta_n $ are close to
$\hat\vartheta,\tilde\vartheta$ random variables and $$\hat\zeta=\sum_{l=1}^{k}
\zeta _l\;{\rm I}\left(\vartheta_l\right)^{-1/2}\;1_{\left\{\HH_l\right\}}.$$

 The proof is based on the weak convergence of
the vector of processes
$$
{\bf Z}_n\left({\bf u}\right)=\left(Z_n^{\left(1\right)}\left(u_1\right),
\ldots,Z_n^{\left(k\right)}\left(u_k\right)\right) ,\qquad
Z_n^{\left(l\right)}\left(u_l\right)=\frac{ L\left(\vartheta
_l+\frac{u_l}{\sqrt{n}},X^T\right)}{L\left(\vartheta
_l,X^T\right)} 
$$
to the limit process ${\bf Z}\left({\bf u}\right)=\left(Z^{\left(1\right)}\left(u_1\right),
\ldots,Z^{\left(k\right)}\left(u_k\right)\right)$, where
$$
Z^{\left(l\right)}\left(u_l\right) =\exp\left\{u_l\Delta _l\left(\vartheta
_l\right)-\frac{u_l^2}{2}  {\rm 
I}\left(\vartheta _l\right)  \right\},\qquad l=1,\ldots,k
$$
(see details in \cite{Kut98}, Section 4.2).

{\bf Example.} Let $\vartheta \in \left(0,3\right)$ and the intensity function
$$
\lambda \left(\vartheta ,t\right)=\left(\vartheta ^3-3\vartheta ^2+2\vartheta
\right) \,t+\left(2\vartheta -3\right)\,t^2+1,\quad 0\leq t\leq 1
$$
then $\lambda \left(1,t\right)=t^2+1$ and $\lambda
\left(2,t\right)=t^2+1$. Hence we have\\ 
$\hat\vartheta _n\Rightarrow
\hat\vartheta=1_{\left\{\left|\zeta _1\right|>\left|\zeta _2\right|\right\}}+
2_{\left\{\left|\zeta _1\right|\leq \left|\zeta _2\right|\right\}}  $ and so on.

\section{Null Fisher information}

Suppose that ${\rm I}\left(\vartheta _0\right)=0$. This means that at one point
$\vartheta _0$ (true value) the function $\dot \lambda\left(\vartheta _0,t\right)=0$
for all $t\in \left[0, \tau \right]$. Moreover,  suppose
that the function $\lambda \left(\vartheta ,t\right)$ is 4 times continuously 
differentiable w.r.t. $\vartheta $ with $\ddot \lambda\left(\vartheta
_0,t\right)=0$ and 
$$
{\rm I}_3\left(\vartheta _0\right)=\int_{0}^{\tau
}\frac{\dddot\lambda\left(\vartheta_0 ,t\right)^2 }{\left(3!\right)^2\;\lambda
\left(\vartheta _0,t\right) }\;
{\rm d}t >0.
$$
Introduce random variable 
$\zeta\left(\vartheta _0\right)\sim {\cal N}\left(0,{\rm I}_3\left(\vartheta
_0\right)\right)$. Then we have:
\begin{align*}
&  n^{1/6}\left(\hat\vartheta _n-\vartheta _0\right)\Longrightarrow \hat u=
\left(\frac{\zeta\left(\vartheta _0\right)}{{\rm I}_3\left(\vartheta
_0\right)} \right)^{1/3}.
\end{align*}

The proof is based on the weak convergence
$$
Z_n\left(u\right)=\frac{L\left(\vartheta_0+
\frac{u}{n^{1/6}},X^n\right)}{L\left(\vartheta_0,X^n\right)}\Longrightarrow
Z\left(u\right)=\exp\left\{u^3\zeta\left(\vartheta _0\right)
-\frac{u^{6}}{2}{\rm I}_3\left(\vartheta _0\right)\right\}.
$$
We have to check the conditions of the Theorem \ref{T2}. Particularly the
estimates \eqref{E}  are replaced by
the estimates
\begin{align*}
c\left|\vartheta _2-\vartheta _1\right|^6\leq \int_{0}^{\tau
}\left[\sqrt{\lambda \left(\vartheta 
_2,t\right)} -\sqrt{\lambda \left(\vartheta _1,t\right)}\right]^{2}{\rm
d}t\leq C\left|\vartheta _2-\vartheta _1\right|^6
\end{align*}

 The limit expression for the bayesian estimator is more
complicated. 

{\bf Example.} Let
$$
\lambda \left(\vartheta ,t\right)=\vartheta \sin^2\left(\vartheta
t\right)+2,\quad 0\leq t\leq 1 ,\quad \vartheta \in \left(-1,1\right)
$$
then ${\rm I}_l\left(0\right)=0,l=1,2 $	 and  ${\rm
I}_3\left(0\right)=\frac{1}{10}. $	Hence 
$$
n^{1/6}\left(\hat\vartheta _n-0\right)\Longrightarrow \left(10\right)^{1/6}
\zeta^{1/3},\qquad \zeta \sim {\cal N}\left(0,1\right)
.$$

\section{Discontinuous Fisher information}

Suppose that the function $\lambda \left(\vartheta ,t\right)$ has at the point
$\vartheta _0$ two different derivatives from the left $\dot 
\lambda \left(\vartheta^{-}_0,t\right)$ and from the right $\dot \lambda
\left(\vartheta^{+}_0,t\right)$ such 
that ${\rm I}\left(\vartheta^{-}_0\right)\neq {\rm
I}\left(\vartheta^{+}_0\right)$ and all the other conditions of regularity are
fulfilled. Then the MLE is consistent, but it is no more asymptotically
normal. Let us introduce a Gaussian vector $\zeta =\left(\zeta _{-},\zeta
_{+}\right)$ with mean zero, $\Ex \zeta _{-}^2=\Ex \zeta _{+}^2=1$ and the covariance 
$$
\Ex\left(\zeta _{-}\zeta _{+}\right) =\Bigl( {\rm I}\left(\vartheta^{-}_0\right){\rm 
                 I}\left(\vartheta^{+}_0\right)  \Bigr)^{-1/2}\;\int_{0}^{\tau
                 }\frac{\dot \lambda \left(\vartheta^{-}_0,t\right)\dot 
                   \lambda \left(\vartheta^{+}_0,t\right) }{\lambda
                 \left(\vartheta_0,t      \right)}{\rm d}t.
$$

Then with the help of the Theorem \ref{T2} it can be shown that the MLE 
 is consistent, and 
$\sqrt{n}\left(\hat\vartheta _n-\vartheta _0\right)\Rightarrow \hat\zeta  $
but its 
limit distribution is a mixture of three random variables:
$$
\hat\zeta=\left\{\begin{array}{ll} 
\frac{\zeta _{-}}{{\rm I}\left(\vartheta^{-}_0\right)^{1/2}}\;, &\;\; {\rm if
}\;\zeta _{-}<0, \zeta _{+}<0$ or $\zeta _{-}<0,\zeta _+>0\; {\rm and
}\;\left|\zeta  _{-}\right|>\left|\zeta_{+} \right|\\ 
0, &\;\; {\rm if}\; \zeta _{-}>0, \zeta _{+}<0,\\
\frac{\zeta _{+}}{{\rm I}\left(\vartheta^{+}_0\right)^{ 1/2}}\;, &\;\;
{\rm if }\;\zeta
_{-}>0, \zeta _{+}>0 $ or $\zeta _{-}<0,\zeta _+>0\; {\rm and }\;\left|\zeta
_{-}\right|<\left|\zeta_{+} \right|
 \end{array} \right.
$$

These properties follow from the form of the limit likelihood ratio process
$$
Z\left(u\right)=\left\{\begin{array}{ll} \exp\left\{u\;\zeta _{-}\;{\rm
I}\left(\vartheta^{-}_0\right)^{1/2}-\frac{u^2}{2}\;{\rm
I}\left(\vartheta^{-}_0\right) \right\},\qquad &\quad u\leq 0\\ 
\exp\left\{u\;\zeta _{+}\;{\rm
I}\left(\vartheta^{+}_0\right)^{1/2}-\frac{u^2}{2}\;{\rm
I}\left(\vartheta^{+}_0\right) \right\},\qquad &\quad u> 0
		       \end{array} \right.
$$
We see that there is an atom at the point 0. This form of the limit likelihood
ratio $Z\left(\cdot \right)$ provides as well the limit distribution of the
bayesian estimates
$$
\sqrt{n}\left(\tilde\vartheta _n-\vartheta _0\right)\Longrightarrow \tilde
u=\frac{\int_{\cal R}^{}uZ\left(u\right){\rm d}u}{\int_{\cal
R}^{}Z\left(u\right){\rm d}u}. 
$$

{\bf Example.} Suppose that $\vartheta \in \left(0,2\right)$ and
$$
\lambda \left(\vartheta ,t\right)=\left(\vartheta
-1\right)\left[3t\,1_{\left\{\vartheta <1\right\}}+ 5 t^2\,1_{\left\{\vartheta
\geq 1\right\}}\right] +15,\quad  0\leq t\leq 1,
$$
 then  ${\rm I}\left(1-\right)=\frac{1}{5}$ and ${\rm
 I}\left(1+\right)=\frac{1}{3}$ and the MLE has the mentioned above limit
 distribution.

\section{Border of the parameter set}

If the true value $\vartheta _0$ is on the border of the parameter set $\Theta
=\left[\alpha ,\beta \right]$, say, $\vartheta _0=\alpha $, then the MLE is
consistent, but 
$$
\sqrt{n}\left(\hat\vartheta _n-\alpha \right)\Longrightarrow 
\frac{\zeta\left(\alpha \right)}{{\rm
I} \left(\alpha\right) } \;1_{\left\{\zeta \geq 0\right\}},\qquad
\zeta\left(\alpha \right) \sim 
{\cal N}\left(0, {\rm I} \left(\alpha\right)\right).
$$
Of course, here ${\rm I} \left(\alpha\right)={\rm I} \left(\alpha^+\right)$.
The estimator is asymptotically half-normal with an atom at 0, i.e., with
probability 0,5 it takes the value 0. This follows from  the form of the limit
likelihood ratio: $$ Z\left(u\right)=\exp\left\{u 
\,\zeta \left(\alpha\right)-\frac{u^2}{2} \, {\rm
I}\left(\alpha\right)\right\},\qquad u\geq 0 .
$$ For the BE we have the limit 
 \begin{align*} 
\sqrt{n}\left(\tilde\vartheta _n-\alpha \right)\Longrightarrow
 \tilde u&=\frac{\int_{0}^{\infty }uZ\left(u\right){\rm d}u}{\int_{0}^{\infty
 }Z\left(u\right){\rm d}u}\\
&= \frac{1}{\sqrt{{\rm I}\,\left(\alpha\right)
 }}\left(\zeta _*+\left(\int_{-\zeta _*}^{\infty }
 e^{-\frac{1}{2}\left(u^2-\zeta _*^2\right)}\,{\rm d}u\right)^{-1}\right).
\end{align*}
where $\zeta _*\sim {\cal N}\left(0,1\right)$.
 
To prove these results we have to check the conditions of the  Theorem
\ref{T2} for the likelihood ratio process
$$
Z_n\left(u\right)=\frac{L\left(\frac{u}{\sqrt{n}},X^n\right)}{L\left(0,X^n\right)},\qquad
u\in \UU_n=\left[0,\beta \sqrt{n}\right]
$$
 with the corresponding limit process.

\section{Cusp type singularity}

Let us suppose that the observed process has intensity function
$$
\lambda \left(\vartheta ,t\right)=a\left|t-\vartheta  \right|^\kappa +\lambda
_0 ,\quad 0\leq t\leq T 
$$
where $\kappa \in \left(0,\frac{1}{2}\right)$. Then this function  is
not differentiable at one point $t=\vartheta $ and the Fisher information 
${\rm I}\,\left(\vartheta \right)=\infty  $. To describe the properties of the
MLE and BE we introduce the normalized  likelihood ratio process
$$
Z_n\left(u\right)=\frac{L\left(\vartheta +\frac{u}{n^{1/2H}
},X^n\right)}{L\left(\vartheta ,X^n\right)},\qquad u\in \UU_n=\left(n^{1/2H}
\left(\alpha -\vartheta _0\right),n^{1/2H}
\left(\beta  -\vartheta _0\right)\right)
$$
and the  limit process 
$$
Z\left(u\right)=\exp\left\{\Gamma_\vartheta W^H\left(u\right)-\frac{\left|u\right|^{2H}}{2}
\Gamma^2_\vartheta\right\},\quad\qquad  u\in {\cal R}.
$$
Here $W^H\left(\cdot \right)$ is double sided fractional Brownian motion, $H=\kappa
+\frac{1}{2}$ (Hurst parameter) and 
$$
\Gamma^2_\vartheta =\frac{4a^2 \sin^2\left(2\pi \kappa \right){\rm B}
\left(1+\kappa,1+\kappa \right)}{ \lambda _0\cos\left(\pi \kappa \right)}\; . 
$$ 
where ${\rm B} \left(1+\kappa,1+\kappa \right)$ is beta function.
 
We can check the conditions of the Theorem 2 and to show that the MLE and BE
are consistent, have the following limits
$$
n^{\frac{1}{2H}} \left(\hat\vartheta _n-\vartheta \right)\Longrightarrow
\hat u, \qquad \quad n^{\frac{1}{2H}}
\left(\tilde\vartheta _n-\vartheta \right)\Longrightarrow \tilde u,
$$
where the random variables are defined by the same equations \eqref{def} and
we have the corresponding convergence of moments.  (For the proof see Dachian, \cite{D}).

\section{Discontinuous  intensity function}

Let us suppose that the observed process
$X^n=\left(X_1\left(\cdot \right),\ldots,X_n\left(\cdot \right)\right)$, where
$X_j\left(\cdot \right)=\left\{X_j\left(t\right),0\leq t\leq T\right\}$ has the
intensity function 
$
\lambda \left(t+\vartheta \right),\quad 0\leq t\leq T
$
and the function $\lambda \left(y\right)$ is positive and  continuously
differentiable everywhere except at the point $\tau $, that is $\lambda
\left(\tau_ +\right)-\lambda \left(\tau_ -\right)=r\neq 0$.  The set $\Theta
=\left(\alpha ,\beta \right)\subset \left(\tau-T ,\tau \right)$. The
likelihood ratio process \eqref{LR} has discontinuous realizations and the MLE
$\hat\vartheta _n$ is defined now by the following relation
$$
\max\left[L\left(\hat\vartheta _n+,X^n\right),L\left(\hat\vartheta
_n-,X^n\right)\right]=\sup_{\vartheta \in \Theta }L\left(\vartheta ,X^n\right).
$$
The BE is defined as before.

  The limit process $Z\left(u\right)$ for the normalized likelihood ratio
$$
Z_n\left(u\right)=\frac{L\left(\vartheta +\frac{u}{n},X^n\right)}{L\left(\vartheta
,X^n\right)},\quad \quad \UU_n=\left(n\left(\alpha -\vartheta
_0\right),n\left(\beta  -\vartheta _0\right)\right) 
$$
is 
$$
Z\left(u\right)=\left\{\begin{array}{ll} \exp\left\{\ln\frac{\lambda
\left(\tau_ +\right)}{\lambda \left(\tau_ -\right)}\,\;\pi
_+\left(u\right)\;-\;\left[\lambda \left(\tau_+ \right)-\lambda
\left(\tau_-\right) \right]\,u \right\},\quad &\quad u\geq 0\\  
\exp\left\{\ln\frac{\lambda
\left(\tau_ -\right)}{\lambda \left(\tau_ +\right)}\,\pi _-
\left(-u\right)-\left[\lambda \left(\tau_+\right)-
\lambda \left(\tau_-\right) \right]\,u \right\},\quad &\quad u\leq 0
		       \end{array} \right.
$$
where $\pi _+ \left(\cdot \right)$ and $\pi _- \left(\cdot \right)$ are
independent Poisson processes of the intensity functions $\lambda \left(\tau_
-\right)$ and $\lambda \left(\tau_ +\right)$ respectively.  Let us denote by
$\hat u$ and $\tilde u $ the random variables defined by the equations
\begin{equation*}
\max\left[Z\left(\hat u+\right),Z\left(\hat
u-\right)\right]=\sup_{u\in {\cal R}}Z\left(u\right),\qquad \quad \tilde u
=\frac{\int_{\cal R}^{}u\,Z\left(u\right)\,{\rm
d}u}{\int_{\cal R}^{}Z\left(u\right)\,{\rm d}u} .
\end{equation*}

 Then the MLE and BE are consistent, have the following limits
\begin{equation*}
\label{12}
n\left(\hat\vartheta _n-\vartheta \right)\Longrightarrow \hat u
,\qquad n\left(\tilde\vartheta _n-\vartheta \right)\Longrightarrow
\tilde u
\end{equation*}
and the convergence of all moments take place. 
It is shown that for all estimators we have a lower bound on the risks 
and the  bayesian estimators are asymptotically efficient. For the proof  see
\cite{Kut98}, section 5.1. 

\section{Windows.} 
 
{\bf Optimal windows}. Suppose that we can have observations on some set
$\BB\subset \left[0,\tau 
\right]$ of Lebesgue measure $\mu\left(\BB\right)\leq \mu _*<\tau $ only. The
family of such sets we denote as ${\cal F}_{\mu _*}$. Our goal is to find the
best window $\BB^*$ and estimator $\vartheta _n^*=\vartheta _n^*\left(
\BB^*\right)$ constructed by the observations $X^n$ on this set $\BB^*$, i.e.;
$X_j=\left\{X_j\left(t\right), t\in \BB^*\right\}$. The best is understood as
the minimizing the mean square error asymptotically
$$
\inf_{\BB\in {\cal F}_{\mu _*}}\inf_{\bar \vartheta _n}\Ex_\vartheta
\left(\bar \vartheta _n\left(\BB \right)-\vartheta \right)^2 \sim \Ex_\vartheta
\left(\vartheta _n^*\left(\BB^* \right)-\vartheta \right)^2.
$$
If we fix the set $\BB$, then  we know that the MLE is asymptotically normal
$$
\sqrt{n}\left(\hat\vartheta _n\left(\BB\right)-\vartheta
\right)\Longrightarrow {\cal N}\left(0, {\rm I}_\BB\left(\vartheta
\right)^{-1}\right),\qquad  {\rm I}_\BB\left(\vartheta
\right)=\int_{\BB}^{}\frac{\dot \lambda \left(\vartheta ,t\right)^2}{\lambda
\left(\vartheta ,t\right)}\;{\rm d}t .
$$
Therefore if we use the MLE then the best $\BB^*=\BB^*\left(\vartheta \right)$
corresponds to the solution of the following equation
$$
{\rm I}_{\BB^*}\left(\vartheta\right)=\sup_{\BB\in {\cal F}_{\mu _*}}{\rm
I}_{\BB}\left(\vartheta\right) .
$$
To solve this equation we introduce the level sets $\CC_{\left(\vartheta
,r\right)} $ and function $\mu 
\left(\vartheta ,r\right)$ as
$$
\CC_{\left(\vartheta ,r\right)}=\left\{t:\;\;\frac{\dot \lambda
\left(\vartheta ,t\right)^2}{\lambda \left(\vartheta ,t\right)}\geq r \right\},\qquad \mu 
\left(\vartheta ,r\right) =\mu \left(\CC_{\left(\vartheta
,r\right)} \right).
$$
Then we define $r_*=r_*\left(\vartheta \right)$ as solution of the equation
$\mu \left(\vartheta ,r \right)=\mu _* $. Now we put
$\BB^*=\CC_{\left(\vartheta ,r_*\right)}$. Of course, we are not obliged to
use the MLE and moreover, this set  $\BB^*$ can not be used for construction
of estimator because it depends on $\vartheta $. Nevertheless it allows to
introduce the lower bound on the risks of all couples ({\it set, estimator}):
for any $\vartheta _0\in \Theta $
$$
\Liminf_{\delta \rightarrow 0 }\Liminf_{n\rightarrow \infty }\;\;\inf_{\BB\in {\cal
F}_{\mu _*},\bar \vartheta _n}\;\;\sup_{\left|\vartheta -\vartheta
_0\right|<\delta }n\,\Ex_\vartheta 
\left(\bar \vartheta _n\left(\BB \right)-\vartheta \right)^2 \geq {\rm
I}_{\BB^*}\left(\vartheta_0\right)^{-1} . 
$$
To construct the asymptotically efficient in this sense couple we first some
 $\BB\in {\cal F}_{\mu _*}$ and by observations
 $X^{\sqrt{n}}=\left\{X_j\left(\BB\right),j=1,\ldots,
 X_{\left[\sqrt{n}\right]}\right\}$ we (consistently) estimate $\vartheta$ using
 some estimator $\bar\vartheta _{\sqrt{n}}$. Then we introduce the observation
 window
$$
\BB_n^*=\left\{t:\;\; \frac{\dot \lambda
\left(\bar\vartheta _{\sqrt{n}} ,t\right)^2}{\lambda \left(\bar\vartheta
_{\sqrt{n}} ,t\right)}\geq r\left(\bar\vartheta _{\sqrt{n}},\mu
_*\right)\right\}. 
$$ 
Now we construct the MLE $\hat\vartheta _{n-\sqrt{n}}$ and show
that 
$$
\lim_{\delta \rightarrow 0 }\lim_{n\rightarrow \infty }\;\;
\sup_{\left|\vartheta -\vartheta
_0\right|<\delta }n\,\Ex_\vartheta  
\left(\bar \vartheta _n\left(\BB_n^* \right)-\vartheta \right)^2 = {\rm
I}_{\BB^*}\left(\vartheta_0\right)^{-1} . 
$$

For the conditions  and proofs see Kutoyants and Spokoiny \cite{KS} or
\cite{Kut98}, Section 4.3.

{\bf Example. } Let 
$
\lambda \left(\vartheta ,t\right)=\left[b+\vartheta
\sin\left(\omega t\right)\right]^2, \; 0\leq t\leq \tau ,
$ 
where $\tau =2\pi /\omega $. Then the Fisher information is
$$
{\rm I}_{\BB}\left(\vartheta\right)=4\,\int_{\BB}^{}\left[\sin\left(\omega
t\right)\right]^2 {\rm d}t .
$$
Introduce 
$$
\CC_r=\left\{t:\;\;  4\sin\left(\omega t\right)^2\geq {r}
\right\},\qquad  r\left(\vartheta ,\mu_* \right)=r_* =4\,\sin^2\left(\frac{2\pi
-\mu_*\omega}{4} \right) 
$$
where  $\left(\mu _*<\tau \right)$.  Then $\varphi
=\arcsin\left(\frac{\sqrt{r_*}}{2}\right)$ 
\begin{align*}
\BB^*&=\left[ \frac{\varphi }{\omega },\;\;  \frac{\tau }{2}-\frac{\varphi }{\omega
}\right]\;\; 
\cup \;\; \left[\frac{\tau }{2}+ \frac{\varphi }{\omega },\;\; \tau-\frac{\varphi
}{\omega }\right]\\
& =\left[\frac{\tau -\mu _*}{4},\;\; \frac{\tau +\mu _*}{4} \right]\cup
\left[\frac{3\tau -\mu _*}{4},\;\; \frac{3\tau +\mu _*}{4} \right] .
\end{align*}
Therefore the observations in the  optimal window are 
$$
X^n=\left(X_1\left(\BB^*\right),\ldots,X_1\left(\BB^*\right) \right),\qquad
{\rm with}\qquad X_j\left(\BB^*\right)=\left\{X_j\left(t\right), t\in \BB^*\right\}
$$
and the asymptotically efficient  estimator is the MLE $\hat\vartheta
_n\left(\BB^*\right)$.  

\bigskip

{\bf Sufficient windows.} The analysis of the proofs of the consistency of the
MLE and BE in the discontinuous case (see, e.g., \cite{Kut98}, Chapter 5)
shows that the main contribution to the likelihood ratio process is made by
the observations near the jumps. Suppose that the intensity function is
$\lambda \left(\vartheta ,t\right)=\lambda \left(t-\vartheta \right)$, where
$\vartheta \in \left(\alpha ,\beta \right)$ and $0<\alpha <\beta <\tau
$. Suppose as well that the function $\lambda \left(s\right), s\in
\left(-\beta ,\tau -\alpha \right)$ is discontinuous 
at some point $\tau _*$ and continuous on $\left[-\beta ,\tau _*)\cup (\tau
_*,\tau -\alpha \right]$. Then it is sufficient to keep the observations on the interval
$\BB=\left[\alpha +\tau _*,\beta +\tau _*\right]$ only, i.e.; to use
$X_j\left(\BB\right)=\left\{X_j\left(t\right), \alpha +\tau _*\leq t\leq \beta
+\tau _*\right\}$, $j=1,\ldots, n$ and the properties of the MLE and BE
(consistency, limit distributions and convergence of moments) will be the same
as in the case of complete observations on $\left[0,\tau \right]$. 

Moreover, if we have a consistent and asymptotically normal estimator $\bar\vartheta _n$ of
$\vartheta $ (say, an estimator of the method of moments), then we can use the
first $\left[\sqrt{n}\right]$ observations for preliminary estimation by
$\bar\vartheta _{\sqrt{n}}$ of the 
window as $\BB_n=\left[\bar\vartheta _{\sqrt{n}}-n^{-1/8},\bar\vartheta
_{\sqrt{n}}+n^{-1/8}\right]$, and then to construct the MLE $\hat \vartheta
_{n-\sqrt{n}}$ and bayesian estimator $\tilde \vartheta
_{n-\sqrt{n}} $. Note that $n^{1/4}\left(\bar\vartheta _{\sqrt{n}}-\vartheta
\right)\Rightarrow  {\cal N}\left(0,\sigma ^2\right)$. Hence 
$$
\Pb_\vartheta \left\{ \left|\bar\vartheta
_{\sqrt{n}}-\vartheta \right| >n^{-1/8}\right\}=\Pb_{\vartheta
}\left\{n^{1/4}\left|\bar\vartheta _{\sqrt{n}}-\vartheta 
\right|>n^{ 1/8}\right\}\longrightarrow 0.
$$ 
Therefore we can have consistent and asymptotically efficient estimators
constructed by observations in the window of vanishing size. In regular case
such effect is difficult to wait. 

{\bf Example.} Let $\vartheta \in \left(\alpha ,\beta \right)\subset
\left(0,\tau \right)$ and
$$
\lambda \left(\vartheta ,t\right)=2at+b\; 1_{\left\{t>\vartheta \right\}},\qquad
0\leq t\leq \tau  .
$$
Then 
$$
\bar\vartheta _{\sqrt{n}}=\tau -\frac{1}{b}\left[\hat\Lambda
_{\sqrt{n}}\left(\tau \right)-a\tau^2 \right] 
$$
is consistent and asymptotically normal estimator of $\vartheta $. Then we
maximize the function
\begin{align*}
L\left(\vartheta ,X^n\right)=&\exp\left\{\sum_{j=\left[\sqrt{n}\right]+1}^{n}
\int_{\bar\vartheta _{\sqrt{n}}-n^{-1/8}}^{\bar\vartheta
_{\sqrt{n}}+n^{-1/8}}\ln \lambda \left(\vartheta ,t\right){\rm
d}X_j\left(t\right)\right.\\
&\quad\left.  -\left(n-\left[\sqrt{n}\right]\right)\int_{\bar\vartheta
_{\sqrt{n}}-n^{-1/8}}^{\bar\vartheta 
_{\sqrt{n}}+n^{-1/8}}\left[ \lambda \left(\vartheta ,t\right)-1\right]{\rm d}t\right\}
\end{align*}
and construct the MLE. Note that the random variable $\bar\vartheta
_{\sqrt{n}} $ is independent on $X_j,j=\left[\sqrt{n}\right]+1,\ldots,n$.

\section{Rates of convergence.} 

 It is interesting to note that if we observe a periodic Poisson process
$X^n=\left\{X\left(t\right), 0\leq t\leq n\right\}$ with the intensity
functions $ \lambda \left(\vartheta+ t\right) $ or $ \lambda \left(\vartheta
t\right) $, where $\lambda \left(\cdot \right)$ is periodic smooth function
(phase and frequency modulations in the optical telecommunication theory),
then we have $\left(n\rightarrow \infty \right)$
$$
\Ex_\vartheta \left(\hat\vartheta _n-\vartheta \right)^2 \sim
\frac{C}{n},\qquad \Ex_\vartheta \left(\hat\vartheta _n-\vartheta \right)^2
\sim \frac{C}{n^3} 
$$
 respectively. If  $\lambda \left(t \right)$ is discontinuous function then
 for the mentioned two cases of modulations we have the different rates
$$
\Ex_\vartheta \left(\hat\vartheta _n-\vartheta \right)^2 \sim
\frac{C}{n^2},\qquad \Ex_\vartheta \left(\hat\vartheta _n-\vartheta \right)^2
\sim \frac{C}{n^4}. 
$$
For the proofs see \cite{Kut98}.

Therefore it is natural to put the following question: {\sl what is the
maximal possible rate of convergence of the mean square error to zero?}
Suppose that we can choose any function $\lambda \left(\vartheta ,t\right),
\vartheta \in \left[0,1\right], t\geq 0$ satisfying the only condition
$$
0\leq \lambda \left(\vartheta ,t\right)\leq L_*,
$$
where $L_*>0$ is some given constant. We denote the class of such functions as
${\cal F}\left(L_*\right)$. It can be shown that
$$
\inf_{\lambda \left(\cdot \right)\in{\cal
F}\left(L_*\right)}\inf_{\bar\vartheta _n}\sup_{\vartheta \in
\left[0,1\right]} \Ex_{\vartheta,\lambda } \left|\bar\vartheta _n-\vartheta
\right|^2=e^{-\frac{nL_*}{6}\left(1+o\left(1\right)\right)},
$$
i.e., the best rate is exponential. To prove this equality we need to prove
two results. The first one is the lower bound for all $\lambda \left(\cdot
\right)\in {\cal F}\left(L_*\right)$ and all estimators $\bar\vartheta _n$
$$
\sup_{\vartheta \in
\left[0,1\right]} \Ex_{\vartheta,\lambda } \left|\bar\vartheta _n-\vartheta
\right|^2\geq e^{-\frac{nL_*}{6}\left(1+o\left(1\right)\right)},
$$
and the second is to construct an intensity function $\lambda _*\left(\cdot
\right)\in {\cal F}\left(L_*\right)$ and an estimator $\vartheta ^*_n$ such
that
$$
\sup_{\vartheta \in
\left[0,1\right]} \Ex_{\vartheta,\lambda _* } \left|\vartheta ^* _n-\vartheta
\right|^2= e^{-\frac{nL_*}{6}\left(1+o\left(1\right)\right)}.
$$
For the proof see Burnashev and Kutoyants \cite{BuK}.

\end{document}